\documentclass[12pt]{article}
\newcommand{\lbl}[1]{\label{#1}}

\begin{document}
\setlength{\baselineskip}{16pt}

\newcommand\tfrac{\frac}
\newcommand\bes{\begin{eqnarray}}
\newcommand\ees{\end{eqnarray}}
\newcommand\bess{\begin{eqnarray*}}
\newcommand\eess{\end{eqnarray*}}
\def \I {\Rightarrow}
\def \N {{\rm  I\!\!\!N}}
\def \R {{\rm  I\!\!\!R}}
\def \Z {{\rm  I\!\!\!Z}}

\newcommand\cL{{\mathcal L}}
\newcommand\cE{{\mathcal E}}
\newcommand\cD{{\mathcal D}}
\newcommand\cM{{\gamma_0}}

\newcommand\bx{{\mathbf x}}
\newcommand\by{{\mathbf y}}
\newcommand\bs{{\mathbf s}}
\newcommand\bn{{\mathbf n}}
\newcommand\rz{{\hbox{\rm z}}}
\newcommand\bI{{\mathbf I}}
\newcommand\bX{{\mathbf X}}
\newcommand\bY{{\mathbf Y}}
\newcommand\bZ{{\mathbf Z}}
\newcommand\rd{{\hbox{\rm d}}}
\def \e {\varepsilon}
\def \ue {u^{\varepsilon}}
\def \ge {g^{\varepsilon}}
\def \Re {R^{\varepsilon}}
\def \uer {u^{\varepsilon}_r}
\def \uet {u^{\varepsilon}_t}
\def \uKe {u^{\varepsilon}_K}
\def \RKe {R^{\varepsilon}_K}
\def \dKe {d^{\varepsilon,K}}
\def \dke {d^{\varepsilon,K}}
\def \uke  {\uKe}
\def \pKe {p^{\varepsilon,K}}
\def \pke  {\pKe}
\def \ske {s^{\varepsilon, K}}
\def \snke {s_n^{\varepsilon,K}}
\def \qke {Q^{\varepsilon,K}}
\def \rKe {r^{\varepsilon,K}}
\def \xKe {x^{\varepsilon,K}}
\def \AKe {A^{\varepsilon,K}}
\def \BKe {B^{\varepsilon,K}}
\def \CKe {C^{\varepsilon, K}}
\def \xe {x^{\varepsilon}}
\def \he {h^{\varepsilon}}
\def \ie {I_{\varepsilon}}
\def \lKe {\lambda^{\varepsilon}_K}
\def \lke {\lambda^{\varepsilon}_K}
\def \a {\lambda}
\def \iom {\int_\Omega}
\def \iqT {\int \int_{Q_T}}
\def \iqt {\int \int_{Q_t}}
\newcommand{\mint}{\displaystyle {\int \kern -0.961em -}}
\newcommand{\mintO}{\displaystyle {\int \kern -0.961em -}_{\Omega}}
\newcommand{\mintG}{\displaystyle {\int \kern -0.961em -}_{\Gamma}}
\newtheorem{thm}{Theorem}
\newtheorem{lemme}{Lemma}
\newtheorem{cor}{Corollary}
\newtheorem{rem}{Remark}
\newtheorem{pro}{Proposition}
\title{Mass conserved Allen-Cahn equation and volume
preserving mean curvature flow}
\author{Xinfu Chen \thanks{Department of Mathematics, University of Pittsburgh,
Pittsburgh, PA 15260, USA.  Part of this work was done during a
visit to the University of Paris-Sud. The author thanks the support
of the National Science Foundation Grant DMS--9971043.}, D. Hilhorst
\thanks{Laboratoire de Math\'ematiques, Analyse Num\'erique et EDP,
CNRS (UMR 8628) et Universit\'e de Paris-Sud, 91405 Orsay Cedex,
France}, E. Logak \thanks{CNRS (UMR 8088) and Department of
Mathematics, Universit\'e de Cergy-Pontoise, 2 rue A. Chauvin, 95302
Cergy-Pontoise Cedex,France}}
\maketitle
\noindent
\begin{abstract}We consider a mass conserved Allen-Cahn equation $u_t=\Delta u+
\e^{-2} (f(u)-\e\lambda(t))$ in a bounded domain with no flux
boundary condition, where $\e\lambda(t)$ is the average of
$f(u(\cdot,t))$ and $-f$ is the derivative of  a double equal well
potential. Given a smooth hypersurface $\gamma_0$ contained in the
domain,  we show that  the solution $u^\e$ with appropriate initial
data  approaches, as $\e\searrow0$, to a limit which takes only two
values, with the jump occurring at the hypersurface
obtained from the volume preserving mean curvature flow
starting from $\gamma_0$.
\end{abstract}

 \maketitle

\section{Introduction.}
\setcounter{equation}{0}
In this paper, we study the limit, as $\e \rightarrow 0$, of the
solution $\ue$ to the mass conserved Allen-Cahn equation $(P^{\e})$
\begin{equation}
 (P^{\e})\left\{
\begin{array}{ll}
u^\e_t=\Delta u^\e +\e^{-2}\big(f(u^\e)- \mintO{f(u^\e)}\big)&\mbox{in } \Omega\times \R^+,\\
\medskip
\partial_\nu u^\e=0 &\mbox{on }
\partial\Omega\times \R^+,\\
\medskip
u^\e(\cdot,0)=\ge (\cdot)&\mbox{on }\Omega\times\{0\},
\end{array}
\right. \lbl{1.1}
\end{equation}
where
$$\mintO{f(\ue)} = \frac{1}{|\Omega|} \int_{\Omega}f(\ue(x,t))dx.$$
Here $\Omega$ is a smooth bounded domain in  $\R^n$ ($n \geq
1$),  $\partial_\nu$  the outward normal derivative to
$\partial\Omega$, and  $-f(u)$ is the derivative of a smooth double
equal well potential;  more precisely,
\begin{equation}\lbl{fpro}   f\in C^\infty(\R),\,
  f(\pm1)=0,\, f'(\pm1)<0, \,\,
\int_{-1}^u f =\int_1^u f <0 \;\;\;\forall u\in(-1,1).
\end{equation}
A typical example is $f(u) = u-u^3$. The initial data $g^\e$
satisfies,  for some smooth hypersurface
$\gamma_0\subset\subset\Omega$,
\begin{equation}\lbl{ge}
\lim_{\e \rightarrow 0} \ge (x) = {\displaystyle{ \left\{
\begin{array}{ll}
-1 & \mbox{ inside } \gamma_0\\
+1 & \mbox{ outside } \gamma_0
\end{array}
\right.\ }}  \quad\forall x\in\bar\Omega\setminus\gamma_0.
\end{equation}
Problem (\ref{1.1}) was proposed, along with its well--posedness,
by Rubinstein and Sternberg \cite{rs} as a
model for phase separation in binary mixture.  The model is mass
preserving and energy decreasing since 
$$\forall t\geq0, \,\,\, \frac{d}{dt}\int_{\Omega} \ue (x,t) dx =0$$
and
$$  \forall t\geq0,  \,\,\, \frac{d}{dt}\iom
\Big(\frac{\e|\nabla \ue|^2}{2} + \frac{1}{\e} F(\ue)\Big)\,dx  =
- \e \int_{\Omega} (u^\e_t)^2  \leq 0,
  $$
 where $\displaystyle{F(u):=-\int_{-1}^u f(s) ds}$ is the double equal well potential.

\noindent
Formally, one can show that, as $\e\to 0$,   the
solution $u^\e$ to (\ref{1.1}) and (\ref{ge}) tends to a limit
\bes\begin{array}{l} \displaystyle{\lim_{\e \rightarrow 0} \ue
(x,t)} = {\displaystyle{ \left\{
\begin{array}{ll}
-1 & \mbox{ inside } \gamma_t\\
+1 & \mbox{ outside } \gamma_t
\end{array}
\right.\ }} \quad\forall x\in\bar\Omega\setminus\gamma_t
\end{array}
\lbl{1.li}\ees where
 $\displaystyle{\Gamma:=\bigcup_{t\geq 0}(\gamma_t\times\{t\})}$
is the solution to the volume preserving mean curvature motion
equation
\begin{equation}
V=(n-1) K_{\gamma_t} - \displaystyle{{(n-1) \over
|\gamma_t|} \int_{\gamma_t}K_{\gamma_t} dH^{n-1}}
 \lbl{1.motion} \mbox{ on } \gamma_t
\end{equation}
starting from  $\gamma_0$. Here $V$ is the normal
velocity of
 $\gamma_t$ (positive when $\gamma_t$ is shrinking) and $K_{\gamma_t}$
 the mean curvature  (positive at points where $\gamma_t$ is locally the boundary of a convex domain).

\noindent
The local in time existence of a unique smooth solution to
(\ref{1.motion}) has been first established in a two-dimensional setting in
\cite{eg}. The general result in arbitrary space dimension is obtained in \cite{escher}, where the
large time behaviour of solutions  for initial data  close to a
sphere was also investigated.  When the initial data is convex,
it is shown in \cite{huisken} that (\ref{1.motion}) admits a
unique global in time convex solution. Related properties of other volume-preserving curvature driven flows are established in \cite{ei}.

\noindent
Concerning the connection between (\ref{1.1}) and (\ref{1.motion}),
Bronsard and Stoth \cite{bs} considered  a radially symmetric case
with multiple interfaces (rings) and proved (\ref{1.li}). Let us
also mention \cite{golo} where a similar result is established for a
different nonlocal mass conserved Allen-Cahn equation, using the
method introduced in \cite{bss}. In the present paper, we shall
consider general smooth initial interfaces $\gamma_0\subset\subset\Omega$ and prove the following:

\medskip

\begin{thm}\lbl{th1}
Let $\Gamma=\bigcup_{0\leq t\leq T}( \gamma_t \times\{t\})$ be a
smooth solution to (\ref{1.motion}) satisfying
$\gamma_t\subset\subset \Omega$ for all $t\in [0, T]$. Then there
exists a family of continuous functions $\{ g^\e\}_{0<\e\leq 1}$
such that the solution $\ue$ to (\ref{1.1})  satisfies (\ref{1.li})
for all $t\in[0,T]$.
\end{thm}

\medskip
\noindent
For the Allen-Cahn equation $u^\e_t=\Delta
u^\e-\e^{-2}f(u^\e)$,  (\ref{1.li}) holds with $\Gamma$ being the
solution to the motion by mean curvature flow $V =
(n-1) K_{\gamma_t}$. A simple method to verify this is to use a
comparison principle and construct sub-super solutions
\cite{c1,ESS}. There are different notions of weak solutions such as viscosity
 \cite{ESS} and varifold \cite{Il} which can be used
to establish the global in time limit.
Nevertheless, (\ref{1.1})
does not have a comparison principle (due to the volume preserving
property) and the simple method does not seem to work. Here we shall
employ a method first used by de Mottoni and Schatzman \cite{ms2}
for the Allen-Cahn equation, and later on by Alikakos, Bates, Chen
\cite{af} for the Cahn--Hillard equation and Caginalp and Chen
\cite{cc} for the phase field system.

\noindent
Namely we first rewrite the equation for $\ue$ in Problem $(P^{\e})$ as
\begin{equation}
u^\e_t=\Delta u^\e +\e^{-2}(f(u^\e)- \e \lambda_{\epsilon}(t))
\mbox{   in } \Omega\times \R^+, \label{eql}
\end{equation}
where we define
\begin{equation}
\forall t \geq 0,\,\,\,\lambda_{\epsilon}(t)=\frac{1}{\e}
\mintO f(u^\e(.,t)) . \label{average}
\end{equation}
The basic strategy of the proof goes as follows.
\begin{enumerate}
\item  For a large enough  $k \in \N$,
construct  an approximate solution $(u^\e_k,\lambda^\e_k)$
satisfying \bes \lbl{1.app} \left\{\begin{array}{ll}
\ue_{k,t}-\Delta \ue_k - \e^{-2} (f(\ue_k)-\e\lambda^\e_k)=
\delta^\e_k \quad  \mbox{ in }  \Omega_T:=\Omega\times [0,T],
\medskip\\     
 \int_\Omega u^\e_{k,t}dx=0 \quad
\forall t\in[0,T], \qquad  \partial_{\nu} \ue_k=0 \mbox{ on }
\partial\Omega_T:=\partial\Omega\times [0,T]
\end{array}\right.
\ees where $\delta^\e_k=O(1)\e^k$. Note that, by integration,
$$\e
\lambda^\e_k=\mintO f(u^\e_k) +O(1)\e^{k+2}. $$

\item For each $t\in[0,T]$ and small positive $\e$, estimate
the lower bound of the spectrum of  the self-adjoint operator $ -
\Delta - \e^{-2} f'(\ue_k(\cdot,t))$; namely,  show that for some
positive constant $C^*$,
\begin{equation}\lbl{spectrap}
\inf_{0<t \leq T} \inf_{0<\e \leq 1}\;\; \inf_{\scriptstyle
\int_\Omega \phi=0, \int_\Omega\phi^2=1} \;\;\iom (|\nabla
\phi|^2 - \e^{-2} f'(\ue_k(.,t)) \phi^2) \geq - C^* .
\end{equation}
\item
Set $R=\ue - u^\e_k$ and show that $R$ tends to $0$ as $\e
\rightarrow 0$.
\end{enumerate}

\noindent
The organization of this paper is as follows. In section $2$, we present an  error
estimate required in step 3.
In section $3$, we recall a known spectrum
estimate \cite{ms,chen} that can be adapted here to prove step $2$ in the strategy described above. After some geometrical preliminary computations in section $4$, we finally construct the approximate solutions in section $5$.

\section{Error Estimate}
The error estimate relies on the following result which is proved in the appendix.
\begin{lemme}\label{gns2}
Let  $\Omega \subset \R^n$ (with $n \geq 1$) be a bounded domain, let $p=\min\{4/n,1\}$. Then there exists $C=C_n(\Omega)>0$ such that  for every
$R\in H^1(\Omega)$ with $\int_\Omega R\,dx=0$,
\begin{equation}\label{interpol}
 \|R\|_{L^{2+p}}^{2+p} \leq C \|R\|_{L^2}^p \|\nabla R\|_{L^2}^{2},
 \end{equation}
where for any $q \geq 1$, $L^q=L^q(\Omega)$.
 \end{lemme}
Rubinstein-Sternberg established in \cite{rs} $L^\infty$ bounds for the solution $\ue$ to Problem
$(P^\e)$ using invariant rectangles. Therefore we can modify $f$ outside of a compact
interval and
assume for simplicity that
$$\lim_{u\to\pm\infty}f(u)=\mp \infty$$
 and that there exists $M>0$ such that
$$ \forall |u| \geq M,\, \, u f''(u) \leq 0. $$
Since $p\in(0,1]$, for any
$C_0>0$, there exists $C=C(C_0,p)$ such that for all $|u| \leq C_0$
and $R\in\R$, 
\bess (f(u+R)-f(u) -f'(u)R) R \leq C |R|^{p+2}.
\eess
Indeed, note that for $R$ in a compact interval,
$$(f(u+R)-f(u) -f'(u)R)R = {f''(u+\theta R) \over 2} R^3 \leq C |R|^{p+2},$$
whereas for $|R|\rightarrow + \infty$,
$f(u +R) R \rightarrow - \infty$, uniformly in  $|u| \leq C_0$ so that
$$(f(u+R)-f(u) -f'(u)R) R \leq (-f(u) -f'(u)R)R \leq C R^2 \leq C |R|^{p+2}.$$
\begin{lemme}\label{ee}
Assume that   $k>\max\{4,n\}$ and
 $\{u_k^{\e}\}_{0<\e\leq
1}$  satisfies   (\ref{1.app}) and
(\ref{spectrap})     with \bess 
 \|\delta_k^\e\|_{L^2(\Omega_T)}
\leq  \e^k, \qquad \|u^\e_k\|_{L^\infty(\Omega_T)}\leq 2. \eess
Let $\{u^\e\}_{0<\e\leq1}$ be solutions to (\ref{1.1}) with
initial data  $\{\ge\}$ satisfying \bess
  \ge(\cdot) = u^\e_k(\cdot,0) + \phi^\e(\cdot),\quad \textstyle
\int_{\Omega} \phi^\e = 0,\quad \|\phi^\e\|_{L^2(\Omega)}\leq
\e^k. \eess Then for all sufficiently small positive $\e$,
\bess
\sup_{0\leq t\leq T} \|\ue(\cdot,t)
-u^\e_k(\cdot,t)\|_{L^2(\Omega)} \leq C(T) \e^k. \eess
\end{lemme}
\begin{rem}
By a bootstrap argument, one can show that other norms of $(\ue
-u^\e_k)$ tend to $0$ as $\e\searrow0$.
\end{rem}
{\it Proof.} In the sequel, $C$ denotes a generic positive
constant independent of $\e$.

\noindent
Set $p=\min\{1,4/n\}$ and $R= \ue - u^\e_k$. Then $\int_\Omega
R(x,t)dx=0$ for all $t\in [0,T]$.   Also,
$$R \{f(\ue)-f(u^\e_k)-f'(u^\e_k) R\} \leq C |R|^{2+p}.$$
 Multiplying by $R$ the
difference of the equations for $u^\e$ and $u^\e_k$ and
integrating the resulting equation over $\Omega$ gives, after
integration by parts, \bess
 &&\frac{1}{2}
\frac{d}{dt}\|R\|_{L^2}^2 +\iom\Big\{ |\nabla R|^2- \e^{-2}
 f'(\uKe)|R|^2\Big\} \leq  \int_\Omega \Big\{C\e^{-2}|R|^{2+p}+|R\,\delta^\e_k|\Big\}.
 \eess
 By (\ref{spectrap}),
  \bess \int_\Omega\Big\{ |\nabla R|^2-\e^{-2} f'(u^\e_k) R^2\Big\}
= \e^2 \int_\Omega \, +(1- \e^2) \int_\Omega \, \geq &\e^2\|\nabla
R\|^2_{L^2} -C \|R\|^2_{L^2}. \eess The interpolation
(\ref{interpol}) then yields
 \bes \frac12\frac{d}{dt} \|R\|^2_{L^2} \leq C \|
\delta^\e_k\|_{L^2} \|R\|_{L^2} + C  \|R\|_{L^2}^2 - \|\nabla
R\|^2_{L^2} \{ \e^2 - C_1 \e^{-2} \|R\|_{L^2}^p \}. \lbl{2.fin}\ees
We define \bess T^\e:=\sup\{ t\in[0,T]  \;|\; \|
R(\cdot,\tau)\|_{L^2} \leq \e^{4/p}C_1^{-1/p}
 \hbox{ \ for all \ }\tau\in [0,t]\}.\eess
Since $k>\max\{4,n\}=4/p$, it follows that
$$\|R(\cdot,0)\|_{L^2}\leq \e^{k}< \e^{4/p} C_1^{-1/p}$$
for $\e>0$ small enough. Therefore, $T^\e>0$. Also, from
(\ref{2.fin}), we have for all $t\in(0,T^\e]$,
 $$ \frac {d}{dt} \|R\|_{L^2}
 \leq C( \|R\|_{L^2} + \|\delta^\e_k\|_{L^2})$$
The Gronwall's inequality then provides
$$\sup_{0\leq t\leq T^\e} \|R(\cdot,t)\|_{L^2}
\leq  e^{CT} [\|R(\cdot,0)\|_{L^2}  +  \int_0^T
\|\delta^\e_k\|_{L^2} dt] \leq C(T) \e^{k}< \tfrac12 \e^{4/p} C_1^{-1/p}$$
 if $\e$ is small enough. Thus, we must have $T^\e=T$. This
completes the proof.

\section{The linearized operator}
\subsection{A Spectrum Estimate}
Assume that $f$ satisfies (\ref{fpro}). Then there is a unique
solution  $\theta_0(\cdot):\R\to (0,1)$  to
\begin{equation}\label{tweq}
\theta_0''+f(\theta_0)=0\mbox{ on } \R,\quad \theta_0(\pm
\infty)= \pm 1,\quad \theta_0(0)=0.\overline{}
\end{equation} The solution
  satisfies, for   $\alpha = \min\{\sqrt{-f'(1)},\sqrt{-f'(-1)}\}$,
$$ D^m_{\rho}\{\theta_0(\rho) \mp 1) =O( e^{-\alpha |\rho|})\hbox{\ \ as \ }\pm \rho \to\infty,
\quad\forall m \in {\N}. $$ Let  $\theta_1 \in C^1(\R)\cap
L^{\infty}(\R)$ be any function satisfying
\begin{equation}\lbl{deftheta1}
 \int_{\R} \theta_0'{}^2\; f''(\theta_0)\;\theta_1\;=0.
\end{equation}
Let $\Omega^-\subset\subset \Omega$ be a subset with $C^3$
boundary $\gamma = \partial \Omega^-$. Denote by $d(x)$  the signed distance
(negative in $\Omega^-$) from $x$ to $\gamma$ and by $s(x)$, for
$x$ close to $\gamma$, the projection from $x$ on $\gamma$
along the normal to $\gamma$.

\noindent
We look for the spectrum of the linearized operator of $-\Delta
u-\e^{-2} f'(u)$ around $u=\psi^\e$ given by
 \bes \psi^\e(x)= \left\{ \begin{array}{ll}
       \theta_0\displaystyle{(\frac{d(x)}{\e}}) + \e\, p^\e(s(x))\theta_1(\displaystyle{\frac{d(x)}{\e}})
+ O(1) \e^2 \quad &\hbox{if \ } |d(x)|\leq \sqrt\e,
\medskip
\\ \pm 1 + O(1)\e &\hbox{if } \pm d(x) \geq \sqrt\e.
\end{array}\right. \lbl{psi}\ees
The following spectrum estimate was first proven by de Mottoni and
Schatzman \cite{ms}, then by Chen \cite{chen} in a more general
situation that can be used in \cite{af,cc}.

\begin{pro}\label{pp}
Let $\gamma \in C^3$, and $p^\e$ and $O(1)$ in (\ref{psi}) be
bounded independently of $\e$. Then there exists a positive constant
$C^*$ depending on $\|\gamma\|_{C^3}$, $\|p^\e\|_{L^\infty}$ and $\|O(1)\|_{L^\infty}$
such that for every $\e \in (0,1]$ and $\phi \in H^1(\Omega)$,
\bess
\iom \Big\{|\nabla \phi|^2 - \e^{-2} f'(\psi^\e) \phi^2\Big\}
\geq - C^* \iom \phi^2 . \eess
\end{pro}
We define the linearized operator around $\theta_0(\rho)$ acting on $v=v(\rho)$ by
\begin{equation}
\cL\; v:= - v''-f'(\theta_0)v \label{lin}
\end{equation}
In our application, $\theta_1$ is the unique solution to
$$\begin{array}{c} \cL\; \theta_1=1- \sigma\,
\theta'_0 \mbox{ \ in \ } \R  \medskip\\
\theta_1(0)=0, \quad\qquad \sigma:= 2/\int_\R \theta_0'{}^2.
\end{array} \label{pb1}$$
 Integrating $\ \theta_0''\;\cL\,\theta_1 \ $ over $\R$ and
by parts, one can verify that (\ref{deftheta1}) is satisfied; see
\cite{ms,chen,af}.

\noindent
We remark that the distance function $d$ in (\ref{psi}) can
be replaced by a ``quasi-distance'' function $d^\e$ given by
$$ {d^\e}(x) = d(x) - \e h_1(s(x))- \e^2 h_2(s(x))+ O(1)\e^3 $$
where $h_1$ and $h_2$ are smooth functions on $s\in\gamma$.
\subsection{Solvability Condition}
\begin{lemme}\label{solv}
Assume that $f$ satisfies (\ref{fpro}). Let $\theta_0$ be the
solution to (\ref{tweq}),
$\alpha=\min\{\sqrt{-f'(1)},\sqrt{-f'(-1)}\}$ and  $\cL$ be defined
in (\ref{lin}). Assume that a function
$h(\rho,s,t)$ satisfies, as
 $\rho \rightarrow \pm \infty$,
$$D_{\rho}^m D_s^n D_t^l[h(\rho,s,t)-h^{\pm}(t)]
=O(|\rho|^i e^{-\alpha |\rho|}) $$ for some $i\geq0$ and all
$(m,n,l)\in \N^3$ and $(s,t)$ in $U \times [0,T]$. Then \bess \cL Q=
h(\cdot,s,t)\hbox{ \ \ in \ } \R,\qquad Q(0,s,t)=0 \eess has a
unique bounded solution $Q(\rho,s,t)$ if and only if
$$ \forall (s,t) \in U \times [0,T],\,\,\,\int_{\R} h(\rho,s,t)\theta'_0(\rho)
d\rho=0.$$ If the solution exists, then it satisfies,  for all
 $(m,n,l)\in \N^3$ and $(s,t) \in U \times [0,T]$,
$$D_{\rho}^m D_s^n D_t^l[Q( \rho,s,t)+
\tfrac{h^{\pm}(t)}{f'(\pm 1)}] =O(|\rho|^i e^{-\alpha |\rho|})\quad
\hbox{as }\rho\to\pm\infty.
$$
\end{lemme}
{\bf Proof} Since $\cL\, \theta_0'=0$, the ode $\cL\,Q=h$ can be
solved explicitly. We omit the details of the proof; see \cite{ms}.

\section{Differential Geometry: local coordinates}
\subsection{Parametrization around the limit interface}
Let
$\Gamma=\displaystyle{\cup_{t\in[0,T]}\;\gamma_t\times\{t\}}\subset\Omega_T$
be the smooth solution to (\ref{1.motion}) on $[0,T]$ and
$\Omega^\pm(t)$ the two domains separated by $\gamma_t$, with
$\gamma_t = \partial \Omega^-(t)$.
 For each fixed $t$, we use
$\rd(x,t)$ to denote the signed distance from $x$ to $\gamma_t$
(positive in $\Omega^+(t) $). Then $\rd(\cdot,\cdot)$ is smooth in a
tubular neighborhood of the interface.
Locally we choose a parametrization
of $\gamma_t$ by $X_0(s,t)$ with $s \in U \subset \R^{n-1}$ so that
\begin{equation}\label{base1}
\big(\frac{\partial X_0}{\partial s_1},...,\frac{\partial X_0}{\partial s_{n-1}}\big)
\end{equation}
is a basis of the tangent space to $\gamma_t$ at $X_0(s,t)$, for each $s \in U$.
We denote
by $\bn(s,t)$ the unit normal vector to $\gamma_t$, pointing towards
$\Omega^+(t)$ so that
$$\bn(s,t)= \nabla \rd(X_0(s,t),t). $$
Up to a suitable multiplication factor $s_1 \rightarrow \lambda s_1$, we may assume that
\begin{equation}\label{det}
\det\big(\bn(s,t), \frac{\partial X_0}{\partial s_1},...,\frac{\partial X_0}{\partial s_{n-1}}\big) =1
\end{equation}
Next for each fixed $t\in [0,T]$, a local parametrization by coordinates $(s,r) \in U \times(-3\delta,3\delta)$ is obtained by
 \begin{equation}
  x=X_0(s,t)+r \,\bn(s,t) = X(r,s,t),     \label{X}
\end{equation}
which defines a local diffeomorphism from $(-3\delta,3\delta)\times U$ onto the tubular neighborhood of $\gamma_t$,
\begin{equation}
V_{3\delta}^t= \{x \in \Omega, \, |d(x,t)| <  3\delta \}.
\label{V3}
\end{equation}
We denote the inverse by
\begin{equation}
\quad
r=\rd(x,t), \,\, s={\bf S}(x,t)=(S^1(x,t),  S^2(x,t),..., S^{n-1}(x,t)).\lbl{3.sh}
\end{equation}
In particular, since for all fixed $s \in U$, $t\in [0,T]$ and for all $r \in (-3\delta,3\delta)$,
$$  d(X_0(s,t)+r \,\bn(s,t), t) =r,$$
it follows by differentiation with respect to $r$ that for all $r \in (-3\delta,3\delta)$,
$$ \nabla d(X_0(s,t)+r \,\bn(s,t), t).\bn(s,t) =1.$$
Using that
\begin{equation}\label{norme}
|\nabla d(x,t)|=1   \mbox{ for } x   \mbox{  close to }\gamma_t,
\end{equation}
this equality imposes that for all $(r,s) \in  (-3\delta,3\delta)\times U$,
\begin{equation}
\nabla d(X_0(s,t)+r \,\bn(s,t), t)= \bn(s,t)
\label{nconst}
\end{equation}
proving that $\nabla d$ is constant along the normal lines to $\gamma_t$.
Thus the projection from $x$ on $\gamma_t$ is defined by
\begin{equation}
X_0({\bf S}(x,t),t)= x-\rd(x,t) \nabla\rd(x,t). \label{xs}
\end{equation}
It follows also from (\ref{norme})
that for all $i = 1, ..., n$ and for $x \in V_{3\delta}^t$,
\begin{equation}\label{sumo}
 \sum_{j=1}^n \frac{\partial^2 d}{\partial x_i
\partial x_j}(x,t)\frac{\partial d}{\partial x_j}(x,t)=0.
\end{equation}
 Thus the symmetric matrix $D_x^2\rd(x,t)$  has eigenvalues
$\{\kappa_1,\cdots,\kappa_{n-1},0\}$ with unit eigenvectors
$\{\tau_1,\cdots,\tau_{n-1},\nabla d\}$ forming an orthonormal basis
of $\R^n$ for $x \in V_{3\delta}^t$. In particular, for $x \in \gamma_t$, the $\tau_i$ are the principal
directions and the $\kappa_i$ are the principal curvatures of
$\gamma_t$. Note that $\{\tau_1,\cdots,\tau_{n-1}\}$ form a basis of
the tangent hyperplane to $\gamma_t$ at $x=X_0(s,t)$.
By definition, $K$ and $K_{\gamma_t}$ are respectively the sum of principal curvatures and the mean curvature of
$\gamma_t$, given by
 \begin{equation}
K=(n-1)K_{\gamma_t}= \Delta d(X_0(s,t),t)=
\sum_{i=1}^{n-1}\kappa_i(s,t). \label{mc}
\end{equation}
Note that using (\ref{sumo}), for $x \in \gamma_t$, we have that
\begin{eqnarray*}
\nabla\rd \cdot\nabla \Delta\rd &=& \sum_{ij}
\frac{\partial}{\partial x_j}\Big(\frac{\partial d}{\partial x_i}\frac{\partial^2 d}{\partial x_i
\partial x_j}\Big)- \sum_{ij} \Big(\frac{\partial^2 d}{\partial x_i \partial x_j}\Big)^2\\
&=& -\sum_{ij}
\Big(\frac{\partial^2 d}{\partial x_i \partial x_j}\Big)^2=- Trace((D_x^2 d)^2) = -\sum_{i=1}^{n-1}
\kappa_i^2.
\end{eqnarray*}
We denote \bes b(s,t)= -\nabla \rd \cdot \nabla
\Delta\rd|_{X_0(s,t),t} = \sum_{i=1}^{n-1}\kappa_i^2.
\lbl{alpha}\ees
Let $V(s,t)$ be the normal velocity of the interface at the point $X_0(s,t)$ so that
using (\ref{nconst}),
\begin{eqnarray}
V(s,t)  &= &X_{0t}(s,t).{\bn}(s,t) \nonumber \\
&=&X_{0t}(s,t). \nabla d(X_0(s,t)+r \,\bn(s,t), t)=- d_t(X(r,s,t),t) \label{vitesse}
\end{eqnarray}
where the last equality follows from differentiating with respect to $t$ the identity
$$ d(X_{0}(s,t)+r \,\bn(s,t), t)=r.$$
It follows that $d_t(x,t)$ is independent of $r=d(x,t)$ for $|r|$ small enough.
Changing coordinates from $(x,t)$ to $(r,s,t)$, we associate to any function $\phi(x,t)$ the function
\begin{equation}
{\tilde \phi}(r,s,t)= \phi(X_0(s,t)+r \bn(s,t),t)
\label{phitilde}
\end{equation}
or equivalently
$$\phi(x,t)= {\tilde \phi}(d(x,t), {\bf S}(x,t),t). $$
By differentiation we obtain the following formulas
\begin{eqnarray}
\partial_t \phi &=& (-V \partial_r + \partial_t^{\Gamma}){\tilde \phi} \nonumber \\
\nabla \phi &=& (\bn \partial_r + \nabla^{\Gamma}){\tilde \phi}  \nonumber \\
\Delta \phi &=& (\partial_{rr} + \Delta d \partial_r  + \Delta^{\Gamma}){\tilde \phi}  \label{chvar}
\end{eqnarray}
with
\begin{eqnarray}
\partial_t^{\Gamma}{\tilde \phi}&=&
(\partial_t + \sum_{i=1}^{n-1} S_t^i \partial_{s^i}){\tilde \phi} \nonumber \\
\nabla^{\Gamma}{\tilde \phi}&=&
(\sum_{i=1}^{n-1}  \nabla S^i \partial_{s^i}){\tilde \phi} \nonumber  \\
 \Delta^{\Gamma}{\tilde \phi} &= &
(\sum_{i=1}^{n-1}  \Delta S^i \partial_{s^i} +
  \sum_{i,j=1}^{n-1}  \nabla S^i . \nabla S^j  \partial_{s^i s^j} ){\tilde \phi} \label{derivtg}
  \end{eqnarray}
where $\nabla S^i$, $S_t^i$, $\Delta d$, $d_t$ are evaluated at
$x= X(r,s,t)$ and are viewed as functions of $(r,s,t)$.
Note that the mixed derivatives of the form $\partial^2_{r s^j} {\tilde \phi}$
do not appear eventually in (\ref{chvar}) because for all $j=1,2,..., n-1$,
$$\nabla S^j(x,t). \nabla d(x,t)=0 $$
(This follows from differentiating with respect to $r$ the identity
$$\forall r \in (-3 \delta, 3 \delta),\,\, S^j(X_0(s,t) + r {\bf n}(s,t),t) = s^j$$
which holds for all fixed $s \in U$, $t\in [0,T]$ and $j=1,2,..., n-1$.)

\subsection{The stretched variable}
Following the method used in \cite{cce}, we now define the stretched variable $\rho$ by considering a graph over
$\gamma_t$ of the form
\begin{equation}\label{gameps}\tilde{\gamma}^{\epsilon}_t=\{ X(r,s,t)  \, / \, r = \e h_{\epsilon}(s,t),\, s \in U\}
\end{equation}
which is (formally) expected to be a representation of the $0$ level set at time $t$ of the solution $\ue$ of
Problem $(P^{\e})$.

\noindent
The stretched variable  $\rho$ is then defined by
\begin{equation}
 \rho = \rho^{\epsilon}(x,t) = \frac{d(x,t) - \e h_{\epsilon}(S(x,t),t)}{\e}
 \label{ro}
 \end{equation}
which represents the distance from $x$ to $\tilde{\gamma}^{\epsilon}_t$ in the normal direction divided by $\e$.
From now on, we use $(\rho, s,t)$ as independent variables for the inner expansions. The relation between the old and new variables are
\begin{eqnarray}\label{Xhat}
  x&=& {\hat X}(\rho, s,t)=X(\e ( \rho + h_{\epsilon}(s,t)),s,t) \nonumber \\
  &=&X_0(s,t)+\e ( \rho + h_{\epsilon}(s,t))\,\,\bn(s,t)
\end{eqnarray}
\noindent We associate to any function $w(x,t)$ the function
\begin{equation}
{\hat w}(\rho,s,t)\,=\,  w(X_0(s,t)+ \e ( \rho +
h_{\epsilon}(s,t))\bn(s,t),t) \label{what}
\end{equation}
or equivalently
$$w(x,t)= {\hat w}(\frac{d(x,t) - \e h_{\epsilon}(S(x,t),t)}{\e},S(x,t),t). $$
Note that
$${\tilde w}(r,s,t)={\hat w}(\frac{r - \e h_{\epsilon}(s,t)}{\e},s,t). $$
By differentiation we obtain the following formulas
\begin{eqnarray}
\partial_t w &=& (-V \e^{-1} - \partial_t^{\Gamma}h_{\e}) {\hat w}_{\rho}
 + \partial_t^{\Gamma} {\hat w} \nonumber\\
\nabla w &= &(\bn \e^{-1}- \nabla^{\Gamma}h_{\e}){\hat w}_{\rho}  + \nabla^{\Gamma} {\hat w} \nonumber \\
\Delta w&=& (\e^{-2}+ |\nabla^{\Gamma}h_{\e}|^2){\hat w}_{\rho \rho}
+ (\Delta d \e^{-1}
- \Delta^{\Gamma}h_{\e}){\hat w}_{\rho}  \nonumber \\
 &&- 2 \nabla^{\Gamma}h_{\e}. \nabla^{\Gamma} {\hat w}_{\rho}
 + \Delta^{\Gamma} {\hat w},  \label{chvar2}
\end{eqnarray}
where in the above formula for $\Delta w$,
\begin{eqnarray}
\Delta d &=& \Delta d|_{x=X_0(s,t)+ \e ( \rho +
h_{\epsilon}(s,t))\bn(s,t)} \nonumber \\
&\approx& K(s,t)-\e ( \rho + h_{\epsilon}(s,t))
b(s,t) + \sum_{i\geq 2}\e^i b_{i}(s,t)( \rho + h_{\epsilon}(s,t))^i,\label{lapd}
\end{eqnarray}
with $b$ defined in (\ref{alpha}),  $K$ defined in
(\ref{mc}) and for some given functions $(b_{i}(s,t))_{i \geq 2}$ only depending on $\gamma_t$.
Therefore
\begin{eqnarray}
 \e^2(\partial_t w - \Delta w)&=& - {\hat w}_{\rho \rho} -\e(V +\Delta d ){\hat w}_{\rho} \nonumber \\
&+& \e^2[ (\partial_t^{\Gamma} {\hat w} - \Delta^{\Gamma} {\hat w})  - ( \partial_t^{\Gamma}h_{\e} - \Delta^{\Gamma}h_{\e}){\hat w}_{\rho}] \nonumber \\
&-& \e^2 [|\nabla^{\Gamma}h_{\e}|^2{\hat w}_{\rho \rho} -  2 \nabla^{\Gamma}h_{\e}. \nabla^{\Gamma} {\hat w}_{\rho}] \label{heat}
\end{eqnarray}
\paragraph{The Jacobi}
For later purposes, we need to compute the Jacobi of the
transformation ${\hat X}$. In the $(\rho,s)$ coordinates, $dx=\e
J^\e(\rho,s,t) ds d\rho$ where $ds$ is the surface element of
$\gamma_t$ and where $\e J^\e(\rho,s,t)=\partial
\hat{X}(\rho,s,t)/\partial(\rho,s)$ is the Jacobi. We prove below
that
\begin{lemme}\label{jacob}
 For all $\rho \in \R$, $s \in U$ and $t \in [0,T]$,
\begin{equation}
J^\e(\rho,s,t)=\prod_{i=1}^{n-1} [ 1+\e (\rho + h^\e(s,t))\kappa_i(s,t)]. \label{jacobi}
\end{equation}
\end{lemme}
{\bf Proof.}
The equality (\ref{jacobi}) is obtained in two steps. First we consider the function $X=X(r,s,t)$ defined in (\ref{X}), denote its Jacobi by $J=J(r,s,t)$ and prove that for all
$\rho \in \R$, $s \in U$ and $t \in [0,T]$,
\begin{equation}\label{jje}
J^\e(\rho,s,t) = J(\e(\rho + h_{\epsilon}(s,t)),s,t).
\end{equation}
Second we compute $J$ and show that for all $s \in U$, for all $t \in [0,T]$,
\begin{eqnarray}
J(r,s,t) & =  & \prod_{i=1}^{n-1} [ 1+ r \kappa_i(s,t)] \nonumber \\
& =  & 1 + \Delta d(X_0(s,t),t) r + \sum_{i=2}^{n-1} r^i j_{i}(s,t),\label{jacobr}
\end{eqnarray}
for some given functions $j_i$ depending on $\gamma_t$.
Consequently (\ref{jacobi}) follows directly from (\ref{jje}) and (\ref{jacobr}).

\noindent
In order to establish (\ref{jje}), note that by definition (\ref{Xhat}),
$${\hat X}(\rho, s,t)=X(\e(\rho + h_{\epsilon}(s,t)),s,t)$$
so that
$$\frac{\partial{\hat X}}{\partial \rho}=\e \frac{\partial{X}}{\partial r}$$
and  for $i=1,..., n-1$,
$$\frac{\partial{\hat X}}{\partial s_i}=\frac{\partial{X}}{\partial s_i} + \epsilon \frac{\partial{\he}}{\partial s_i}\frac{\partial{X}}{\partial r}. $$
Thus for all
$\rho \in \R$, $s \in U$ and $t \in [0,T]$,
\begin{eqnarray*}
&&\e J^\e(\rho,s,t)=\e \det \big[\frac{\partial{X}}{\partial r},\frac{\partial{X}}{\partial s_1}+ \epsilon \frac{\partial{\he}}{\partial s_1}\frac{\partial{X}}{\partial r},...,\frac{\partial{X}}{\partial s_{n-1}} + \epsilon \frac{\partial{\he}}{\partial s_{n-1}}\frac{\partial{X}}{\partial r}\big] \\
&&= \e \det \big[\frac{\partial{X}}{\partial r},\frac{\partial{X}}{\partial s_1},...,\frac{\partial{X}}{\partial s_{n-1}}\big](\e(\rho + h_{\epsilon}(s,t)),s,t) = \e J(\e(\rho + h_{\epsilon}(s,t)),s,t)
\end{eqnarray*}
which is (\ref{jje}).

\noindent
In order to establish (\ref{jacobr}), we consider the Hessian matrix of $d$ on $\gamma_t$ and
denote for $s \in U$ and $t \in [0,T]$
$$A= A(s,t) = D_x^2\rd(X_0(s,t),t)$$
so that (\ref{sumo}) reads
\begin{equation}\label{sumo2}
 A.{\bf n}(s,t) = 0 .
 \end{equation}
 Moreover, differentiating the identity (\ref{nconst}) at $r=0$ with respect to $s_i$ for $i=1,..., n-1$ yields
\begin{equation}\label{nderiv}
 A.\frac{\partial{X_0}}{\partial s_i} = \frac{\partial{\bn}}{\partial s_i} .
 \end{equation}
From
$$ X(r,s,t)\, =\,X_0(s,t)+r \,\bn(s,t)  , $$
it follows that using (\ref{sumo2})
$$\frac{\partial{X}}{\partial r}=\bn(s,t) = (I_n + r A(s,t))(\bn(s,t))$$
and that, using (\ref{nderiv})  for $i=1,..., n-1$,
$$\frac{\partial{X}}{\partial s_i}=\frac{\partial{X_0}}{\partial s_i} + r \frac{\partial{\bn}}{\partial s_i} =
(I_n + r A(s,t))(\frac{\partial{X_0}}{\partial s_i}). $$
Therefore for all $s \in U$ and $t \in [0,T]$,
\begin{eqnarray*}
&&J(r,s,t)=\det \big[\frac{\partial{X}}{\partial r},\frac{\partial{X}}{\partial s_1},...,\frac{\partial{X}}{\partial s_{n-1}}\big] \\
&&= \det \big[(I_n + r A)(\bn),(I_n + r A)(\frac{\partial{X_0}}{\partial s_1}),...,
(I_n + r A)(\frac{\partial{X_0}}{\partial s_{n-1}})\big]\\
&&= \det \big[I_n + r A(s,t)\big] \det \big[\bn,\frac{\partial{X_0}}{\partial s_1},...,\frac{\partial{X_0}}{\partial s_{n-1}}\big]
\end{eqnarray*}
which in view of (\ref{det}) proves that
$$J(r,s,t)= \det \big[I_n + r A(s,t)\big]$$
which yields (\ref{jacobr}), since the eigenvalues of $A(s,t)$ are
$\kappa_1,\cdots,\kappa_{n-1},0$.

\section{The approximate solution}

\subsection{Asymptotic Expansions}
Let $k>\max\{2,n/2\}$ be a fixed  integer. In the sequel, we use the
sign $\approx$ to represent an asymptotic expansion; namely,
$\displaystyle{\phi^\e\approx \sum_{i\geq0}\e^i \phi_i}$ means that for every
integer $j \in \N$, $\displaystyle{\phi^\e=\sum_{i=0}^{j} \e^i \phi_i +
O(1)\e^{j+1}}$ where $O(1)$ is bounded independently of $\e\in(0,1)$.
For example, since $f$ is smooth, for any
bounded sequence $\{b,a_0, a_1, a_2, ... \}$, we have the asymptotic
expansion \bes f(b + \e \sum_{i\geq0}\e^i a_i)&\approx&
\sum_{j\geq0} {\e^j f^{(j)}(b)} \Big(\sum_{i\geq0} \e^i
a_i\Big)^j /j\,! \nonumber
\\ &\approx&  f(b)+\e f'(b)\sum_{i\geq0} \e^i a_i +
\e^2\sum_{i\geq0} \e^i f_i(b,a_0,..., a_i)    \lbl{fi} \ees
where for any fixed $b$, $f_i(b,a_0,..., a_i)$ is a polynomial in $(a_0,..., a_i)$
of degree $\leq i$.
\paragraph{Outer expansion}
We expand $\lambda^\e(t)$ and $u^\e(x,t)$ for
$|\rd(x,t)|\geq 3\delta$ by
 \begin{equation}\label{lambdae}
 \lambda^\e(t)\approx\lambda_0(t) + \e\,\lambda_1(t)+\e^2\,\lambda_2(t)+\cdots
 \end{equation}
 \begin{equation}\label{uepm}
 \ue(x,t) \approx u^\pm_\e(t):=\pm 1 + \e \{u_0^\pm(t) + \e
u_1^\pm(t) +\cdots\}.
\end{equation}
Substituting (\ref{lambdae}) and (\ref{uepm}) into (\ref{eql}) gives
  $$f(u^\pm_\e(t))=\e \lambda^\e(t) + \e^2 (u^\pm_{\e})'(t)$$
  which yields
 for all $i\geq 0$,
\bes   u_{i}^\pm(t) =\{ \lambda_{i}-f_{i-1}(\pm1,u^\pm_0,\cdots,u^\pm_{i-1})
- u^\pm_{i-2,t}\}/f'(\pm1) \lbl{upm} \ees                   where
$f_{-1}=u^\pm_{-2}=0$, $u^\pm_{-1}=\pm1$, and  $f_i$ ($i\geq0$)
is as in (\ref{fi}).  Hence, $u^\pm_i$ are  determined by
$\{\lambda_0,\cdots,\lambda_i\}$.
\paragraph{Inner expansion}
We shall assume that $h^\e$ has the asymptotic expansion
 \begin{equation}\label{he}
 \e h^\e(s,t) \approx \e h_1(s,t)+ \e^2 h_2(s,t)
+ \cdots,\quad (s,t)\in U \times [0,T]
\end{equation}
Near the interface, we assume that
the  function ${\hat \ue}$ associated to $\ue$ by (\ref{what}) has the asymptotic expansion
 \begin{equation}\label{inner}
  {\hat \ue}(\rho,s,t) \approx \theta_0(\rho) + \e \, \{u_0(\rho,s,t)
 + \e\,u_1(\rho,s,t) + \cdots\}.
 \end{equation}
In the sequel, the zero-th order expansion refers to
$$\{\rd(x,t),\lambda_0(t),u_0(\rho,s,t),u_0^\pm(t)\}$$ and the $i$-th
order expansion refers to 
$$\{h_i(s,t),\lambda_i(t),
u_i(\rho,s,t),u^\pm_i(t)\}.$$
 We shall use $(\cdots)_{i-1}$ to denote a
generic function of $(\rho,s,t)$ depending only on expansions of order
 $\leq i-1$.
\paragraph{Matching condition}
We impose that for all $i \in  \N$,
\begin{equation}
 \forall (s,t)\in U\times[0,T],\,\, \lim_{\rho \rightarrow \pm \infty} u_i(\rho,s,t) = u^\pm_i(t)
 \label{match}
 \end{equation}
\paragraph{Translation}
We also impose for all $i \in \N$,
\begin{equation}
\forall (s,t)\in U\times[0,T],\,\, u_i(0, s,t) =0
\label{trans}
 \end{equation}
to be consistent with the assumption that $\rho=0$ is the $0$-level set of $\ue$.
\subsection{The u-equation in the new variables}
The equation (\ref{eql}) is
$$ -f(u) =  -\e^2 (u_t - \Delta u) -\e  \lambda_{\e}(t). $$
In the new variables $(\rho, s,t)$, using (\ref{heat}), it becomes the following equation for the function $u={\hat \ue}$ associated to $\ue$ by
(\ref{what}),
\begin{eqnarray}
-f(u)&=& u_{\rho \rho} +
\epsilon[( V(s,t)  + \Delta d)  u_{\rho}   - \lambda_{\e}] \nonumber \\
 &+&  \epsilon^2[ (\Delta^{\Gamma}u
 -\partial_t^{\Gamma} u ) + ( \partial_t^{\Gamma}h_\epsilon - \Delta^{\Gamma}h_\epsilon)u_{\rho} ] \nonumber\\
 &+& \epsilon^2[|\nabla^{\Gamma}h_\epsilon|^2 u_{\rho \rho} - 2 \nabla^{\Gamma}h_\epsilon. \nabla^{\Gamma}u_{\rho}] ,\label{new}
\end{eqnarray}
where $V(s,t)$ is given by
(\ref{vitesse}) and $\Delta d$ is expanded from (\ref{lapd}) as
\begin{equation}
\Delta d \approx  K(s,t)  - \sum_{i\geq 1}   \e^i [b(s,t) h_i(s,t) +
 \delta_{i-1}(\rho, s,t)],\label{laplace}
\end{equation}
with $\delta_{i-1}$ depending only on expansions of order $\leq i-1$
 (in particular, $\delta_{0}(\rho, s,t)= \rho b(s,t)$). Note
 that $\delta_{i-1}(\rho,s,t)$ is a polynomial in $\rho$ of degree $\leq i$, whose
 coefficients are polynomial in $(h_1,..., h_{i-1})$ with $(s,t)$-dependent
 coefficients.
 \subsection{The recursive i-th equations}
\paragraph{The zeroth order expansion}
Since $\theta_0$ defined in (\ref{tweq}) satisfies
$$-f(\theta_0)= (\theta_0)_{\rho \rho},  \,\,\, \theta_0(\pm \infty) = \pm 1, \,\,\, \theta_0(0) =0, $$
the equation (\ref{new}) is satisfied at zeroth order as well as the matching and translation condition
(\ref{match})-(\ref{trans}).
\paragraph{The first order expansion}
At first order $(\e^1)$, the equation (\ref{new}) imposes
\begin{equation}\label{u1}
\cL \, u_0= (K(s,t) + V(s,t))(\theta_0)'(\rho) - \lambda_{0}(t),
\end{equation}
with $\cL$ defined in (\ref{lin}).
The solvability condition stated in Lemma \ref{solv} reads
$$(K(s,t) + V(s,t) )\int_{\R}(\theta'_0)^2(z) dz= 2 \lambda_{0}(t)$$
which reads in view of (\ref{pb1})
\begin{equation}
V(s,t) =  - K(s,t) + \sigma \lambda_{0}(t) \mbox{ for } s \in U
\label{solv1}
\end{equation}
which implies in view of (\ref{vitesse}) that
\begin{equation}
d_t =  \Delta d  -\sigma \lambda_{0}(t) \mbox{ on } \gamma_t.
\label{eqd}
\end{equation}
Moreover equation (\ref{u1}) has then a unique solution satisfying
(\ref{match})-(\ref{trans})which is given by
\begin{equation}
u_0(\rho,s,t)=-\lambda_{0}(t)\theta_1(\rho)
\label{u0}
\end{equation}
for all $(s,t) \in U \times [0,T]$. Note that for all non-negative
$m,n,l$,
$$D_{\rho}^m D_s^n D_t^l[u_0(\rho,s,t)-u_0^{\pm}(t)]
=O( e^{-\alpha |\rho|}) \mbox{ as } \rho \rightarrow \pm\infty.$$
\paragraph{Higher order expansion}
Plugging the expansions (\ref{fi}), (\ref{he}), (\ref{inner}) into (\ref{new}) and using (\ref{solv1}) and (\ref{laplace})
leads to the following identity
\begin{eqnarray}
&&- f(\theta_0)-\e f'(\theta_0)(\sum_{i\geq0} \e^i u_i) -
\e^2\sum_{i\geq0} \e^i f_i(\theta_0,u_0,..., u_i)  \nonumber \\
 &=& \theta_0'' + \e (\sum_{i\geq0} \e^i (u_i)_{\rho \rho}) +
\epsilon[(\sigma \lambda_0(t)- \sum_{i\geq 1}   \e^i (b h_i +
 \delta_{i-1}) )u_{\rho}   - \sum_{i\geq0} \e^i \lambda_i] \label{L1} \\
 &+&  \epsilon^3 \sum_{i\geq0} \e^i (\Delta^{\Gamma}
 -\partial_t^{\Gamma})u_i  - \e(\sum_{i\geq1} \e^i (\Delta^{\Gamma}
 -\partial_t^{\Gamma})h_i)(\theta_0' + \e \sum_{i\geq0} \e^i (u_i)_{\rho})  \label{L2} \\
 &+& [\e^2|\nabla^{\Gamma} h_\epsilon|^2 u_{\rho \rho} - 2\epsilon(\sum_{i\geq1} \e^i \nabla^{\Gamma}
h_i). \nabla^{\Gamma}u_{\rho}].  \label{L3}
\label{newi}
\end{eqnarray}
Define the operator ${\cal N}^{\Gamma}$ acting on
functions $h=h(s,t)$ by
\begin{equation}
{\cal N}^{\Gamma} h:=(\partial_t^{\Gamma}h-\Delta^{\Gamma}h- b h) \label{Ng}
\end{equation}
We derive below the $(i+1)$-th order expansion for $i\geq 1$ and obtain the following result.
\begin{lemme}\label{expi}
At order $\e^{i+1}$, with $i\geq 1$, the equation (\ref{new}) imposes
\begin{equation}
\cL \, u_i={\cal N}^{\Gamma}(h_{i}) \theta_0' - \lambda_{i}(t)
 +b_{12} (\nabla^{\Gamma}h_1.\nabla^{\Gamma}h_{i})\theta_0''+ R_{i-1}(\rho,s,t),\label{equi}
\end{equation}
with $R_{i- 1}$ only depending on expansions of order $\leq i-1$.
Besides $R_{i-1}(\rho,s,t)$ is a polynomial in $\rho$ of degree
$\leq i$ (whose
 coefficients are polynomial in $(h_1,..., h_{i-1}, u_1,...,
 u_{i-1})$ and in their derivatives with respect to $(\rho, s,t)$).
\end{lemme}
{\bf Proof.}
First note that using (\ref{solv1}), the coefficient of order $\e^{i+1}$ in (\ref{L1}) is
\begin{eqnarray}
&&(u_i)_{\rho \rho} + \sigma \lambda_0(t)(u_{i-1})_{\rho}  -
b(s,t)h_{i}(s,t)\theta_0' - \lambda_{i}(t) + (\cdots)_{i-1}
\nonumber \\
&&(u_i)_{\rho \rho}  - b h_{i} \theta_0' - \lambda_{i}(t) +
(\cdots)_{i-1}, \label{la}
\end{eqnarray}
with $(...)_{i-1}$ depending only on expansions of order $\leq i-1$. Moreover in view of (\ref{laplace}),
it is a polynomial in $\rho$ of degree $\leq i$ (whose
 coefficients are polynomial in $(h_1,..., h_{i-1}, u_1,...,
 u_{i-1})$ and in their derivatives with respect to $(\rho, s,t)$).

\noindent
Next, in view of (\ref{he}), the coefficient of order
$\e^{i+1}$ in (\ref{L2}) is
\begin{eqnarray}
&&(\Delta^{\Gamma}
 -\partial_t^{\Gamma}) u_{i-2} +   ( \partial_t^{\Gamma} -\Delta^{\Gamma})h_{i}\theta_0'+ (...)_{i-2} \nonumber \\
 &=& ( \partial_t^{\Gamma} -\Delta^{\Gamma})h_{i} \theta_0'+ (\cdots)_{i-2}.
 \label{lb}
 \end{eqnarray}
To obtain the term
of order $\e^{i+1}$ in (\ref{L3}), note that
\begin{eqnarray*}
\e^2|\nabla^{\Gamma}h_\epsilon|^2  &\approx &|\sum_{i\geq 1} \e^{i} \nabla^{\Gamma}h_{i}|^2 \approx \sum_{i\geq 2}
\e^{i}(\sum_{j=1}^{i-1}\nabla^{\Gamma}h_j.\nabla^{\Gamma}h_{i-j}) \\
&\approx &\e^2 |\nabla^{\Gamma}h_1|^2 +  \sum_{i\geq 3} \e^{i} (2 \nabla^{\Gamma}h_1.\nabla^{\Gamma}h_{i-1} +(\cdots)_{i-2})
\end{eqnarray*}
so that
$$\e^2|\nabla^{\Gamma}h_\epsilon|^2 \approx [\e^2 |\nabla^{\Gamma}h_1|^2 +  \sum_{i\geq 3} \e^{i} (2 \nabla^{\Gamma}h_1.\nabla^{\Gamma}h_{i-1} +(\cdots)_{i-2})][\theta_0'' + \e \sum_{i\geq0} \e^i (u_i)_{\rho \rho} ]. $$
Hence the coefficient of order $\e^{i+1}$ in $\e^2|\nabla^{\Gamma}h_\epsilon|^2 u_{\rho \rho}$ is
$$b_{1,2} (\nabla^{\Gamma}h_1.\nabla^{\Gamma}h_{i})\theta_0'' + (\cdots)_{i-2}$$
with $b_{1,2} =1$ or $2$ for $i=1$ or $i\geq 2$ respectively.

\noindent
Similarly, the coefficient of order $\e^{i+1}$ in the term
$\e^2 \nabla^{\Gamma}h_\epsilon. \nabla^{\Gamma}u_{\rho}$ is
$$\nabla^{\Gamma}h_{i-1}. \nabla^{\Gamma}(u_0){\rho} + \nabla^{\Gamma}h_{i-2}. \nabla^{\Gamma}(u_1){\rho} + ... +
\nabla^{\Gamma}h_{1}. \nabla^{\Gamma}(u_{i-2}){\rho}$$
where the first term cancels out since $\nabla^{\Gamma}(u_0)'=0$ in view of (\ref{u0}); thus it only depends on expansions of order $\leq i-2$ and appears below in the remainder.

\noindent
Finally at order $\e^{i+1}$, with $i\geq 1$, the equation (\ref{new}) reads
\begin{eqnarray*}
&&-f'(\theta_0)u_i -f_{i-1}(\theta_0,u_0,..,u_{i-1})=(u_i)_{\rho\rho}  - \lambda_{i}(t)\\
&&
 + (\partial_t^{\Gamma}h_{i}-\Delta^{\Gamma}h_{i} - b h_{i}) \theta_0'
 +b_{12} (\nabla^{\Gamma}h_1.\nabla^{\Gamma}h_{i})\theta_0''+ R_{i-1}(\rho,s,t),
\end{eqnarray*}
with $R_{i- 1}$ only depending on expansions of order $\leq i-1$.
Moreover $R_{i-1}(\rho,s,t)$ is a polynomial in $\rho$ of degree
$\leq i$ as described in Lemma \ref{expi}.
\paragraph{The solvability condition}
According to Lemma \ref{solv}, the equation (\ref{equi}) has a
solution if and only if the following solvability condition is
satisfied.
\begin{equation}\label{ortho}
\forall (s,t) \in U \times [0,T],\,\,\,\int_{\R}\cL \, u_i(\rho,s,t)
\theta_0'(\rho) d\rho =0.
\end{equation}
Note that
$$ \int_{\R}\,  b_{12} (\nabla^{\Gamma}h_1.\nabla^{\Gamma}h_{i})\theta_0''(\rho) \theta_0'(\rho) \, d\rho =
 b_{12} (\nabla^{\Gamma}h_1.\nabla^{\Gamma}h_{i})(s,t)\int_{\R}\, \frac{1}{2}[(\theta_0')^2]'(\rho) \, d\rho =0 $$
so that the condition (\ref{ortho}) reads
\begin{equation}
{\cal N}^{\Gamma}(h_{i}) = \sigma \lambda_{i}(t)  + r_{i-1}(s,t),
\label{solvi}
\end{equation}
with
$$r_{i-1}(s,t) =- \frac{\sigma}{2} \int_{\R}\,R_{i-1}(\rho,s,t) (\theta_0)'(\rho) \, d\rho $$ 
only depending on expansions of order $\leq (i-1)$. We summarize these
results in the next lemma.
 \begin{lemme}\lbl{44}
 Let $k \geq 1$ be given. Assume that for all $i \leq k-1$,
(\ref{equi}) has a solution $u_i$ satisfying
\begin{equation}\lbl{424}
D_{\rho}^m D_s^n D_t^l[u_i(\rho,s,t)-u_i^{\pm}(t)] =O( \rho^i
e^{-\alpha |\rho|}) \mbox{ as } \rho \rightarrow \pm \infty.
\end{equation}
Also assume that for $i=k$,  $\{h_i(s,t),\lambda_i(t)\}$ satisfies
(\ref{solvi}). Then for $i=k$, (\ref{equi}) admits a unique solution
satisfying $u_i(0,s,t)=0$ and (\ref{424}).
\end{lemme}
The proof follows from Lemma \ref{solv} and an induction argument
and is omitted. 
Just note that in the limit $\rho\to\pm\infty$, the
equation $0=\e^2 (u^\e_t-\Delta u^\e)+f(u^\e)-\e
\lambda^\e|_{x={\hat X}(\rho,s,t)}$ 
becomes the outer expansion
equation, so that $u_i(\pm\infty,s,t)=u_i^\pm(t)$. Furthermore since $R_{i-1}$ is a polynomial in
$\rho$ of degree $\leq i$, (\ref{424}) is
satisfied for each $i\geq 0$ and $(s,t)\in U\times [0,T]$.
\subsection{Equation for $\lambda^\e$.}
To find  $\lambda^\e(t)$, we use an asymptotic expansion for
$0=\int_\Omega u^\e_t(x,t)dx$. We denote by $\Omega^\pm_\e(t)$ the
two domains separated by $\tilde{\gamma}^\epsilon_t$ defined in
(\ref{gameps}), with $\tilde{\gamma}^\epsilon_t = \partial
\Omega^-_\e(t)$. Hence in view of (\ref{ro})
\begin{eqnarray}
  \Omega^+_\e(t) &=& \{ x\in\Omega\;|\;d(x,t)> 3 \delta \} \, \cup \, \{ x\in V_{3 \delta}^t\;|\;
  [d(x,t)- \e h^\e({\bf S}(x,t),t]>0\}\nonumber \\
  &=&\{ x\in\Omega\;|\;d(x,t)> 3 \delta \} \, \cup \,\{ x\in
  V_{3 \delta}^t\;|\;\rho^{\e}(x,t)>0\} \label{omegap}
\end{eqnarray}
and
\begin{eqnarray}
  \Omega^-_\e(t) &=& \Omega \setminus {\overline{\Omega^+_\e(t)}} \nonumber \\
  &=&\{ x\in\Omega\;|\;d(x,t)<- 3 \delta \} \, \cup \,\{ x\in
  V_{3 \delta}^t\;|\;\rho^{\e}(x,t)<0\} \label{omegam}
\end{eqnarray}
We write
\begin{equation}
\int_\Omega u^\e_t(x,t)dx = \int_{|d(x,t)|\geq 3\delta} u^\e_t(x,t)dx +
\int_{|d(x,t)|< 3\delta} u^\e_t(x,t)dx \label{utlp}
\end{equation}
where
\begin{equation}
\int_{|d(x,t)|< 3\delta} u^\e_t(x,t)dx=\int_{|\rho^{\e}(x,t) |\geq \frac{\delta}{\e}}  u^\e_t(x,t)dx +\int_{|\rho^{\e}(x,t)|<\frac{\delta}{\e}}  u^\e_t(x,t)dx
\label{utp}
\end{equation}
In the sequel we choose $0 <\e \leq \e_0$ small enough so that
\begin{equation}\label{ehe}
\forall \e \in (0, \e_0], \,\,\, \max_{s \in U, t \in [0,T)} |\e \he (s,t)| \leq \frac{\delta}{2}
\end{equation}
Then it follows that
$$|\rho^{\e}(x,t) |\geq \frac{\delta}{\e}  \,\, \Rightarrow \,\, |d(x,t)|\geq \frac{\delta}{2}.$$
Thus if  ${|d(x,t)| \geq 3\delta}$ or  $|\rho^{\e}(x,t) |\geq \frac{\delta}{\e} $,  then
$|d(x,t)| \geq {\delta \over 2}$ so that at these points $(x,t)$,
 $$u^\e_t(x,t) \approx (u^\e_+)'(t)\chi_{\{d(x,t)>0\}} +(u^\e_-)'(t)\chi_{\{d(x,t)<0\}}$$
  (exponentially small terms of order $O(e^{-\alpha\delta \over 2\e})$ do not affect the asymptotic expansion in the $\e$ power series). Therefore in view of (\ref{utlp})-(\ref{utp})
  \begin{eqnarray}
&& \int_\Omega u^\e_t(x,t)dx \approx
\int_\Omega [ (u^\e_+)'(t)\chi_{\{d(x,t)>0\}} +(u^\e_-)'(t)\chi_{\{d(x,t)<0\}}] \,dx \\
&&+ \int_{|\rho^{\e}(x,t)|<\frac{\delta}{\e}}  
  [ u^\e_t - (u^\e_+)'(t)\chi_{\{d(x,t)>0\}} - (u^\e_-)'(t)\chi_{\{d(x,t)<0\}} ]\,dx 
  \nonumber \\
&&\approx
I_1  
+ \int_{|\rho^{\e}(x,t)|<\frac{\delta}{\e}}  
  [ u^\e_t - (u^\e_+)'(t)\chi_{\{d(x,t)>0\}} - (u^\e_-)'(t)\chi_{\{d(x,t)<0\}} ] dx ,
\label{coupe}
\end{eqnarray}
 where
\begin{equation}
I_1= (u^+_{\e})'(t)|\Omega^+_\e(t)| +
 (u^-_{\e})'(t)|\Omega^-_\e(t)|.\label{I1}
\end{equation}
In the second integral, we make the change of variables given in
(\ref{Xhat}) and substitute the expression of $\ue_t$ in formula (\ref{chvar2}) to obtain
\begin{eqnarray}
 &&\int_{|\rho|<\delta/\e} [\ue_t - (u^\e_+)'(t)
\chi_{\{d(x,t)>0\}} - (u^\e_-)'(t)\chi_{\{d(x,t)<0\}} ] \,dx = \nonumber \\
 && \int_{0<\rho<\delta/\e} {\partial_t^{\Gamma}}[\hat{\ue} (\rho,s,t) - (u^\e_+)(t)] \epsilon
  J^\epsilon(\rho,s,t) \, d\rho\,  ds  \nonumber \\
  &&+
  \int_{-\delta/\e<\rho<0} {\partial_t^{\Gamma}}[\hat{\ue} (\rho,s,t) - (u^\e_-)(t)] \epsilon
  J^\epsilon(\rho,s,t) \, d\rho \, ds \nonumber \\
&&+ \int_{|\rho|<\delta/\e} (-V \e^{-1} -
\partial_t^{\Gamma}h_{\e})\frac{\partial{\hat{\ue}}}{\partial \rho} \epsilon
  J^\epsilon(\rho,s,t) \, d\rho \, ds
   \label{tint}
\end{eqnarray}
Finally,
 $$\int_\Omega u^\e_t(x,t)dx \approx  I_1 +I_2 + I_3,$$ where
\begin{equation}
I_2=  \int_{|\rho|<\delta/\e} {\partial_t^{\Gamma}}[\hat{\ue}
(\rho,s,t) - (u^\e_+)(t)\chi_{\{\rho>0\}} -
(u^\e_-)(t)\chi_{\{\rho<0\}} ] \epsilon
  J^\epsilon(\rho,s,t) \,d\rho \,ds \label{I2}
\end{equation}
and
\begin{equation}
I_3= \int_{|\rho|<\delta/\e} (-V \e^{-1} -
\partial_t^{\Gamma}h_{\e})\frac{\partial{\hat{\ue}}}{\partial \rho} \epsilon
  J^\epsilon(\rho,s,t) \, d\rho \, ds.\label{I3}
\end{equation}
{\it The calculation for $I_1$.} The boundary of
$\Omega_\e^-(t)$ is $\tilde{\gamma}^\epsilon_t$ which according to (\ref{gameps}) is given in local coordinates $(r,s)$ by $r = \e \he(s,t)$. Therefore in view of (\ref{jacobr}), we have that
\begin{eqnarray*}
|\Omega^-_\e(t)| &=& |\Omega^-(t)| + \int_{U}\int_0^{\e\he(s,t)}
J(r,s,t) \, dr \, ds  \\
 &\approx& |\Omega^-(t)| +
\sum_{i\geq 1} \e^i \{ \int_{U}\,  h_i(s,t) \, ds + (...)_{i-1}  \},
\end{eqnarray*}
where $(...)_{i-1}$ only depends on expansions of order $\leq i-1$.
Hence
\begin{eqnarray*}
|\Omega^+_\e(t)|&=&|\Omega|-|\Omega^-_\e(t)| \\
 &\approx& |\Omega^+(t)| -
\sum_{i\geq 1} \e^i \{ \int_{U}\,  h_i(s,t) \, ds + (...)_{i-1}  \}.
\end{eqnarray*}
 From the outer expansion,
$$u^\pm_{\e, t}  \approx \e\sum_{i\geq 0} \e^{i}
(u_{i}^{\pm})'(t)\approx \sum_{i\geq 1} \e^{i}
(u_{i-1}^{\pm})'(t),$$
with $(u_{i-1}^{\pm})'(t) $ given by (\ref{upm}) and depending only on expansions of
order $\leq i-1$.
Therefore
\bess I_1=u^+_{\e,t}(t)|\Omega^+_\e(t)| +
u^-_{\e,t}(t)|\Omega^-_\e(t)|   \approx  \Sigma_{i\geq 1}\e^i
(...)_{i-1} \eess where $(...)_{i-1}$ depends only on expansions of
order $\leq i-1$.

\noindent
 {\it The calculation for $I_2$.}
 Using the expression for ${\partial_t^{\Gamma}}\hat{\ue}$ in formula (\ref{derivtg}) and (\ref{424}), we compute
\begin{eqnarray*}
&&{\partial_t^{\Gamma}}[\hat{\ue}(\rho,s,t) - (u^\e_+)(t)\chi_{\{\rho>0\}} -
(u^\e_-)(t)\chi_{\{\rho<0\}} ] \\
&&\approx \e\sum_{i \geq 1} \e^i {\partial_t^{\Gamma}}[u_i(\rho,s,t) - u_i^+(t)\chi_{\{\rho>0\}}
-  u_i^-(t)\chi_{\{\rho<0\}}] \\
&&\approx \e\sum_{i \geq 1} \e^i (\partial_t + \sum_{j=1}^{n-1} S_t^j \partial_{s^j})[u_i(\rho,s,t) - u_i^+(t)\chi_{\{\rho>0\}}
-  u_i^-(t)\chi_{\{\rho<0\}}]  \\
&&\approx \sum_{i \geq 2} \e^i  O(\rho^{i-1} e ^{- \alpha |\rho|})
\end{eqnarray*}
with  $ O(\rho^{i-1} e ^{- \alpha |\rho|})$ depending only on expansions of
order $\leq i-1$. Therefore by definition of $I_2$ in (\ref{I2}),
 $$I_2 \approx \sum_{i \geq 3} \e^i  (...)_{i-2},$$
where $(...)_{i-2}$ depends only on expansions of
order $\leq i-2$.

 \noindent
 {\it The calculation for $I_3$.}
  Using the expansions
 $$\frac{\partial{\hat{\ue}}}{\partial \rho}
  \approx \theta_0'+\e\sum_{i\geq 0}\e^i \frac{\partial{u_i}}{\partial \rho},$$
  $$(-V  -
\e \partial_t^{\Gamma}h_{\e})=
  d_t(X_0(s,t),t)-\sum_{i\geq 1}\e^i{\partial^\Gamma_t h_i}$$
and rewriting the expression of $ J^\epsilon$ in (\ref{jacobi}) as
\begin{eqnarray*}
J^\e(\rho,s,t)&=&\prod_{i=1}^{n-1} [ 1+\e (\rho + h^\e(s,t))\kappa_i(s,t)] \\
&\approx& 1 + \Delta d(X_0(s,t),t) \e( \rho +h^\e(s,t)) + \sum_{i\geq
2}\e^i (...)_{i-1},
\end{eqnarray*}
with  $(...)_{i-1}$ depending only on expansions of
order $\leq i-1$, we obtain that
\begin{eqnarray*}
&&(-V  -\e \partial_t^{\Gamma}h_{\e}) \frac{\partial{\hat{\ue}}}{\partial \rho}J^\e(\rho,s,t) \approx \\
&&d_t(X_0(s,t),t) \theta_0'(\rho)+\sum_{i\geq 1}\e^i \theta_0'(\rho)(-\partial_t^\Gamma h_i+
d_t(X_0(s,t),t) h_i \Delta\rd )+ \sum_{i\geq 1}\e^i (...)_{i-1}
\end{eqnarray*}
so that
\begin{eqnarray*}
 I_3&\approx &\int_{U}\int_\R  \Big\{
 \theta_0' d_t(s,t) +\Sigma_{i\geq1}\e^i [\theta_0' (-\partial_t^\Gamma h_i+
d_t(s,t)\Delta\rd(s,t) h_i)+(...)_{i-1}]\Big\} d\rho ds
 \\ &\approx & 2 \int_{U} d_t(s,t)\,ds + \sum_{i\geq1} \e^i  \Big\{
2\int_{U} \{-\partial^\Gamma_t h_i +
(\rd^t\Delta\rd)h_i\}\,ds + (...)_{i-1}\Big\}.
\end{eqnarray*}
Finally, substituting $\rd_t$ and $\partial^\Gamma_t h_i$ by
(\ref{eqd}) and (\ref{solvi}), and using
$\int_{U}\Delta^{\Gamma} h_ids=0$,    we obtain
 \bess \tfrac12\int_\Omega u^\e_t\approx  \int_{U}
 (\Delta \rd-\sigma\lambda_0) \,ds
 + \sum_{i\geq1}\e^i \Big\{ \int_{U} [(-b
+\rd_t\Delta\rd)h_i-\sigma\lambda_i]     \,ds + (...)_{i-1}
\Big\}\eess
Thus the condition $\int_\Omega u^\e_t\,dx\approx 0$ is
equivalent to
 \bes  \sigma\lambda_0(t) &=& \overline{\Delta \rd(\cdot,t)},
 \lbl{lambda0}\\
\sigma\lambda_i(t) &=& -\overline{[b(\cdot,t)
-\rd_t(\cdot,t)\Delta\rd(\cdot,t)]h_i(\cdot,t)}+
\Lambda_{i-1}(t),  \quad i\geq1 \lbl{lambda}\ees where
$\Lambda_{i-1}(t)$ depends only on expansions of order $\leq
i-1$, and $ \overline{\phi(\cdot)} := \frac{1}{|U|}
\int_{U} \phi $, the average of $\phi$ over $\gamma_t$ parametrized by $U$.

\noindent
Hence, we obtain closed systems for $\rd$,$h_1$,$\cdots$, $h_i$, namely
 \bes
&&\rd_t(s,t)= \Delta \rd(s,t)-\overline{\Delta
\rd(s,t)}, \label{5.d0}\\ 
&&\partial^\Gamma_t h_i = \Delta^\Gamma h_i
+b \,h_i -\overline{[b (\cdot,t)
-\rd_t(\cdot,t)\Delta(\cdot,t)]h_i(\cdot,t)}+\Lambda_{i-1}(t)       \lbl{5.h0}
 \ees
on $U\times[0,T]$.
\subsection{Construction of Expansions of Each Order}
We can now use induction to construct each order of expansion as
follows:

1) {\it Zeroth order}. Given a smooth initial interface $\gamma_0$,  (\ref{5.d0}) is
equivalent to the volume preserving mean curvature flow
(\ref{1.motion}). By the result established in \cite{escher}, there is a time $T>0$ such that there is a unique smooth
solution on a time interval $[0,T]$. Consequently, $\Gamma=\bigcup_{0\leq t\leq T}(\gamma_t \times\{t\})$
and the modified distance function $d$ are well defined. Set
$\lambda_0(t)$ as in (\ref{lambda0}), $u_0(\rho,s,t)$ as in (\ref{u0}) and  $u_0^\pm(t) = \lambda_0/f'(\pm1)$ as in (\ref{upm}). We obtain
the zeroth order expansion $\{\rd(x,t),\lambda_0(t),u_0(\rho,s,t),u_0^\pm(t)\}$.

2) {\it Higher order expansion}. Fix $i\geq 1$. Assume that  all
expansions of order $\leq i-1$ are constructed.  Then
$\Lambda_{i-1}(\cdot)$ in (\ref{5.h0}) is  known. Since $\gamma_t$  is
a smooth hypersurface without boundary, it follows from standard parabolic PDE theory \cite{LUS} that
(\ref{5.h0}) admits a unique smooth solution (assuming an initial
condition such as $h_i(\cdot,0)=0$ on $U$ is given).
Consequently, we can define $\lambda_i(t)$ as in (\ref{lambda}),
 $u^\pm_i$ as in (\ref{upm}) and $u_i$ as the solution of (\ref{equi}) given by Lemma \ref{44}.
This gives  the $i$-th order expansion $\{h_i(s,t),\lambda_i(t),
u_i(\rho,s,t),u^\pm_i(t)\}$ and  completes the induction.

\subsection{Construction of the Approximate Solution}
Now fix an arbitrary positive integer $k$. We  construct an
approximate solution $u_k^\e$ such that Lemma \ref{ee} can be
applied.

\noindent
Let $\delta>0$ be a small fixed constant such that (i)
$\rd(x,t)$ is smooth in a $3\delta$-neighborhood of $\Gamma$, and
(ii) for each $t\in[0,T]$, $\gamma_t$ is at least $3\delta$ distance
away from $\partial\Omega$. We define \bess
   \rho_k^\e (x,t) &=& \e^{-1}\{ \rd(x,t)-\Sigma_{i=1}^{k+1} \e^i
h_i({\bf S}(x,t),t)\},\\
 u_{\e,k}^{{in}}(x,t) &=& \theta_0(\rho_k^\e) + \e\sum_{i=0}^{k+1} \e^i u_i(\rho^\e_k(x,t),{\bf S}(x,t),t) ,
\\  u_{\e,k,\pm}^{out}(t) &=& \pm 1+ \e \sum_{i=0}^{k+1} \e^i u_i^\pm(t).
\eess
   We note that $\rho_k^\e$,$u_{\e,k}^{in}$ are smooth  in
a $3\delta$ neighborhood of $\Gamma$.

\noindent
Now let $\zeta(s)\in C^\infty(\R)$ be a cut-off function
(depending only on $\delta$) satisfying 
\begin{eqnarray*}
&&\zeta(s)=1\hbox{ \
if \ }|s|\leq \delta,\quad \zeta(s)=0 \hbox{ \ if \
}|s|>2\delta, \\
&&  \quad 0\leq s\zeta'(s) \leq 4 \hbox{ \  if \ }
\delta\leq|s|\leq 2\delta.
\end{eqnarray*}
We define the needed approximation solution $u^\e_k$ by  \bess
\tilde u^\e_k(x,t)&:=& \zeta(\rd) u^{in}_k + [1-\zeta(\rd)]
\Big\{ u^{out}_{\e,k,+}  \chi_{\{\rd>0\}} + u^{out}_{\e,k,-}
\chi_{\{\rd<0\}}\Big\},
\\
u^\e_k(x,t)&:=& \tilde u^\e_k(x,t) + \mintO \{ \tilde
u^\e_k(.,0)-\tilde u^\e_k(.,t) \}
 \eess for all
 $(x,t)\in\bar\Omega\times[0,T]$.  Then by construction 
 $u^\e_k$ is an approximation of order $k$ as needed in Lemma
\ref{ee}. Here we just remark that  (i) in the set $\{(x,t)\;|\;
\delta\leq\pm\rd(x,t)\leq 2\delta\}$, the limiting behavior
(\ref{424}) guarantees that
$u^\e_k(x,t)=u^{out}_{\e,k,\pm}(t)+O(e^{-\alpha\delta/(4\e)})$,
valid also for differentiation, (ii) $\partial_n u^\e_k=0$ on
$\partial\Omega_T$ since $u^\e_k$ is   a function of $t$  near
$\partial\Omega_T$, and (iii) the correction 
$$\int_\Omega \{ \tilde u^\e_k(.,0)-\tilde u^\e_k(.,t)\}= -
 \int_\Omega  \int_{[0,t]}  (\tilde u^\e_k)_t (y,\tau) \, d \tau  \, dy$$ is of order $O(\e^{k+1})$,
valid also for differentiation.

\noindent
This completes the construction of the approximating solution, and
also the proof of Theorem \ref{th1}.

\paragraph{Appendix A}: Proof of Lemma \ref{gns2}.

\noindent
We first consider the case $n \geq 4$ so that $p= 4/n$. The
Gagliardo-Nirenberg-Sobolev inequality (see \cite{evans}, Theorem 2, p.265) states that there exists $C>0$ such that for every
$R\in H^1(\Omega)$,
$$\|R\|_{L^{2^*}} \leq C \|R\|_{H^1},$$
with $ 2^* =\displaystyle{\frac{2n}{n-2}}$. Using Poincar\'e-Wirtinger inequality (see \cite{evans}, Theorem 1, p.275), it follows that there exists $C>0$ such that for every
$R\in H^1(\Omega)$ with $\int_\Omega R\,dx=0$,
\begin{equation}\label{ev}
\|R\|_{L^{2^*}} \leq C \|\nabla R\|_{L^2}.
\end{equation}
Writing H\"older inequality, we have that
$$ \|R\|_{L^{2+p}}^{2+p}= \int_{\Omega} |R|^2 |R|^{p} \leq (\int_{\Omega}|R|^{2\beta})^{1/\beta}(\int_{\Omega}|R|^{p\beta'})^{1/\beta'}$$
and we choose
$$\beta =  \frac{n}{n-2} = \frac{2^*}{2}, \,\,\, \beta' =  \frac{n}{2}$$
to obtain
$$ \|R\|_{L^{2+p}}^{2+p} \leq  \|R\|_{L^{2^*}}^{2} \| R\|_{L^2}^{p}.$$
Combined with (\ref{ev}), this yields the inequality
$$ \|R\|_{L^{2+p}}^{2+p} \leq C \|R\|_{L^2}^p \|\nabla R\|_{L^2}^{2},$$
which is the conclusion of Lemma \ref{gns2}.

\noindent
Next we consider the case that $1 \leq n \leq 3$ so that $p =1$.
Schwarz's inequality then gives that
$$ \|R\|_{L^{3}}^{3}= \int_{\Omega} |R|^2 |R| \leq
\|R\|_{L^{4}}^{2} \|R\|_{L^{2}}$$
For $n=1,2,3$, by Sobolev's imbedding theorem, $H^1 \subset L^4$,  so that
there exists $C>0$ such that for every
$R\in H^1(\Omega)$,
$$\|R\|_{L^{4}} \leq C \|R\|_{H^1}.$$
Using again Poincar\'e-Wirtinger inequality, we finally deduce that there exists $C>0$ such that for every
$R\in H^1(\Omega)$ with $\int_\Omega R\,dx=0$,
$$ \|R\|_{L^{3}}^{3} \leq
C  \|\nabla R\|_{L^2}^{2}  \|R\|_{L^{2}},$$
which concludes the proof of Lemma \ref{gns2}.

\end{document}